# Distributionally Robust Joint Chance-Constrained Optimization Framework for Electricity Imbalance: Integrating Renewables and Storages

A. Noori, B. Tavassoli*, A. Fereidunian

Faculty of Electrical Engineering, K.N. Toosi University of Technology, Tehran, Iran.
E-mail: anoori@email.kntu.ac.ir, tavassoli@kntu.ac.ir, Fereidunian@eetd.kntu.ac.ir
*Corresponding author

**Abstract**
Integrating Distributed Energy Resources (DERs) with peer-to-peer (P2P) energy trading offers promising solutions for grid modernization by incentivizing prosumers to participate in mitigating peak demand. However, this integration also introduces operational uncertainties and computational challenges. This paper aims to address these challenges with a novel scalable and tractable distributionally robust joint chance-constrained (DRJCC) optimization framework that effectively facilitates P2P energy trading by enhancing flexibility provision from large-scale DER operations under uncertain supply and demand. Therefore, a practical framework is proposed to solve the core challenges of DRJCC by integrating three key components: (1) a Wasserstein ambiguity set that effectively quantifies uncertainty with sparse data, (2) a CVaR-based approximation of joint chance constraints to balance computational efficiency with risk control, and (3) a privacy-preserving ADMM algorithm that enables distributed implementation through decomposition. To discern patterns in the data that indicate collaboration potential and adjust ambiguity sets for improved efficiency, K-means clustering is applied to historical scenarios. Simulation results show that the proposed framework reduces peak demand by approximately 28% and total community costs by around 31%, underscoring its effectiveness in enhancing grid robustness, operational reliability, and economic optimization in renewable-based energy management.

NOMENCLATURE
Notations          Description

*Sets and Indices*
$m, n$          Indices for prosumers
$i, j$          Indices for data samples
$k$          Index for iterations
$h$          Index for time periods
$N, N_n$          Set of prosumers, neighbours of $n$
$I$          Size of the dataset
*Parameters*
$p_n^g$          Photovoltaic power output (kW)
$p_n^l$          Power of must-run loads (kW)
$c_p, c_{nm}$          Predicted DA and P2P energy prices ($/kWh)
$\gamma_n^b, \gamma_n^s, \gamma_n^p$          Coefficients of cost functions
$\overline{p}_n^b, \overline{p}_n^s, \overline{p}_n, \overline{q}_n$          The maximum power limits (kW)
$\underline{p}_n^b, \underline{p}_n^s, \underline{p}_n, \underline{q}_n$          The minimum power limits (kW)
$\underline{E}_n, \overline{E}_n$          The limits of state-of-charge (kWh)
$\underline{S}_n, \overline{S}_n$          The limits of energy shift state (kWh)
$\eta$          The coefficient of battery storage system
$\rho$          The radius of the ambiguous set
$\sigma$          Augmented Lagrangian parameter
$\varepsilon$          Tolerance for the ADMM algorithm's residual
*Decision Variables*
$p_n^b$          Power of battery storage system (kW)
$p_n^s$          Consumption power of shiftable loads (kW)
$p_n$          DA energy exchange with the main grid (kW)
$q_n$          RT energy exchange with the main grid (kW)
$p_{nm}$          P2P energy trade of prosumers $n$ and $m$ (kW)

$P_n$          Total P2P energy trade (kW)
$E_n$          SoC of battery storage system (kWh)
$S_n$          Energy shift state of shiftable loads (kWh)

## 1. Introduction

### 1.1. Motivation and Background
The Iranian power grid faces a growing supply-demand gap, particularly during peak summer periods, exacerbated by urbanization, economic growth, and the adoption of electricity-intensive technologies, leading to severe energy deficits, with a shortfall of 14,000-17,000 megawatts in the summer of 2024 [1]. In response, many consumers and prosumers are enhancing their energy resilience by integrating renewable sources like solar PV and investing in energy storage. This presents an opportunity for the government to reduce grid investments while supporting the decentralized energy transition. P2P energy trading integrated with DERs incentivizes prosumers to contribute flexible resources and engage in collective grid-supportive behaviours, but this requires effective coordination in a distributed, scalable framework under uncertainty. By fostering growth in DERs and crafting supportive policies, the government can promote a decentralized, resilient grid that enhances energy security and stability.

## 1.2. Related Works

The successful transition to a sustainable power grid necessitates the seamless integration of renewable energy sources and energy storage systems. While these technologies offer significant potential for decarbonization, their inherent intermittency introduces considerable complication in grid operations, potentially leading to imbalances and jeopardizing system reliability. To address this challenge, a range of approaches have been proposed in the literature, including stochastic programming (SP), robust optimization (RO), and distributionally robust optimization (DRO). These methods exhibit significant variations in their treatment of uncertainty, ranging from scenarios with full knowledge of the underlying probability distributions to those characterized by bounded uncertainty and diverse forms of ambiguity sets.

On the other hand, decentralized and adaptive approaches are essential for effectively integrating renewable energy sources. Decentralization promotes local energy sharing, reducing reliance on centralized grid infrastructure, and enhancing resilience to local energy demands and resource variability [2]. Adaptive strategies, supported by data analytics, enable real-time peak shaving and optimize energy consumption, improving grid stability and facilitating demand-response mechanisms [3]. Moreover, risk-aware approaches are critical for managing the uncertainties inherent in renewable generation [4]. These advancements highlight the need for scalable, uncertainty-driven solutions that combine decentralization, adaptability, and risk management to create a resilient and sustainable power grid.

From a data-driven distributionally robust optimization perspective, these challenges are addressed with varying focuses and considerations. In [6], uncertainty in renewable generation is addressed using a chance-constrained model for joint energy and reserve markets. It applies a quadratic programming (QP) approach and versatile distribution modeling, focusing on peer-to-peer (P2P) energy trading. However, its centralized structure limits scalability for larger networks. As a seminal work on distributionally robust joint chance constraints (JCC) for energy and reserve dispatch [7], employs Wasserstein DRO to ensure reliability but remains centralized, targeting static, day-ahead (DA) planning with affine policies. A data-driven DRO framework for co-optimizing P2P energy trading and network operation in interconnected microgrids is introduced in [8]. A distributed approach with affine policies improves scalability. Extending DRO to local renewable energy aggregators [9] adopts a Weibull distribution model for uncertainty. While comprehensive in scope, its centralized implementation poses limitations in real-time applications.

The work [10] explores ambiguous chance constraints with perturbation-based approximations for microgrid energy management. It adopts a centralized structure for real-time (RT) operation, leveraging Chernoff's inequality (CI) for tractability. The paper [11] uniquely integrates chance-constrained co-optimization for P2P trading and grid operation. By using Gaussian distributions and individual chance constraints (ICC), it provides efficient RT solutions for localized networks. Focusing on integrated electricity and heating systems, [12] adopts Wasserstein DRO with JCCs. A distributed conic programming (CP) approach enhances scalability and decision-making reliability. The paper [13] merges stochastic DRO with Stackelberg game theory for energy hubs, adopting clustering techniques to reduce computational complexity in centralized real-time systems.

These studies focus on integrating DERs using distributionally robust optimization and stochastic chance-constrained approaches, with some addressing P2P energy trading and grid support. However, they lack scalable, privacy-preserving, and holistic frameworks for decentralized systems. Challenges include balancing real-time computational efficiency with robust risk management, reliance on rigid assumptions for uncertainty quantification, and a predominant focus on individual chance-constraints and day-ahead planning. Furthermore, the potential of clustering techniques to improve collaboration and efficiency in P2P trading remains underexplored.

## 1.3. Contributions

In this paper, we propose a novel distributed distributionally robust joint chance-constrained optimization framework that integrates the following key features:

- *Real-time adaptability* for P2P energy trading and flexibility provision in dynamic grid operations,
- *Privacy-preserving distributed implementation* using a consensus-ADMM algorithm,
- *Clustering techniques* to enhance efficiency and collaboration potential among participants.

A detailed comparison of the proposed framework with state-of-the-art methods is provided in Table 1.

Table 1. Comparison of related literature

| Ref. | P2P energy | Grid support | Uncertainty | Optimization | (St- Dy) / (DA- RT) | Constraints | Scalability | ML |
|------|-----------|-------------|-------------|--------------|---------------------|-------------|-------------|-----|
| [5]  | –  | DR | W | MILP | Dy/DA | AP  | Cen. | – |
| [6]  | ✓  | PF | V | QP   | St/DA | ICC | Dis. | – |
| [7]  | –  | R  | W | CP   | St/DA | JCC | Cen. | – |
| [8]  | ✓  | PF | W | LP   | Dy/DA | AP  | Dis. | – |
| [9]  | –  | PF | W | QP   | St/DA | Wei | Cen. | – |
| [10] | –  | PF | P | MILP | Dy/RT | CI  | Cen. | – |
| [11] | ✓  | PF | G | QP   | St/RT | ICC | Cen. | – |
| [12] | –  | PF | W | CP   | Dy/DA | JCC | Dis. | – |
| [13] | –  | R  | W | MILP | Dy/RT | AP  | Cen. | C |
| [14] | –  | –  | W | MILP | Dy/DA | ICC | Dis. | – |
| This paper | ✓ | PD | W | CP | Dy/RT | JCC | Dis. | C |

AP: affine policy, C: clustering, Cent: centralized, CI: Chernoff's inequality, CP: conic program, DA: day-ahead, Dis: Distributed, DR: demand response, Dy: Dynamic, G: Gaussian, IA: inner approx., ICC: individual chance-constrained, JCC: joint chance-constrained, LP: linear program, MILP: mixed-integer linear program, P: perturbation, PD: peak demand, PF: power flow, QP: quadratic program, R: reserve, RT: real-time, St: Static, V: Versatile distribution, Wei: Weibull distribution, W: Wasserstein.

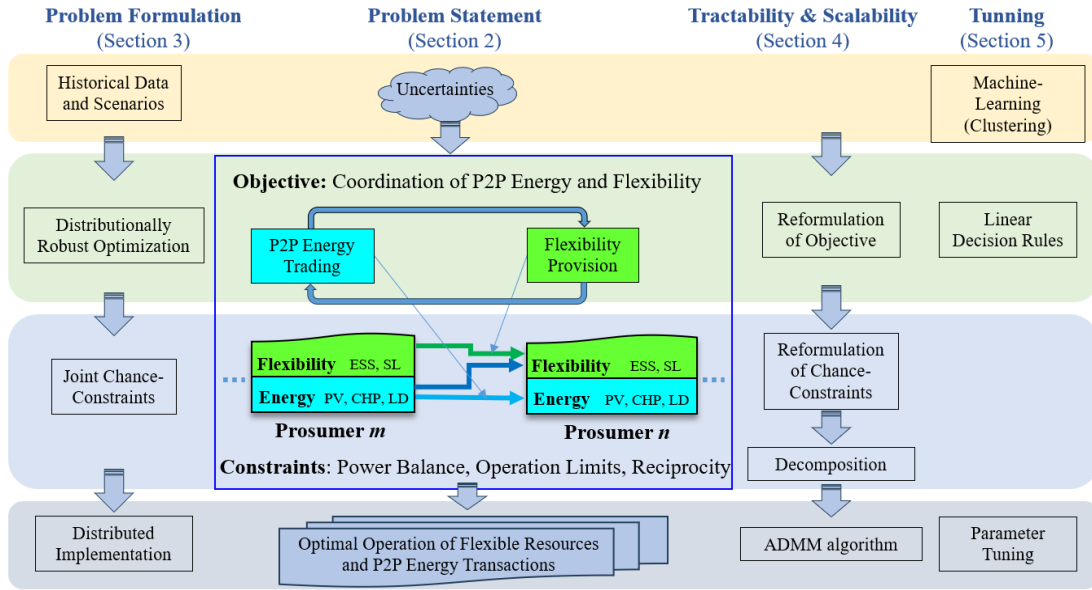

Fig. 1 Proposed distributed distributionally robust joint chance-constrained optimization framework for coordination of P2P energy trading and distributed flexibility provision

The remainder of this paper is organized as follows: Section 2 presents the proposed framework and model. Section 3 formulates the two-stage distributionally robust joint chance-constrained problem. Subsequently, Section 4 focuses on problem reformulation to enhance computational tractability and scalability. Section 5 presents a case study with numerical simulations and Section 6 discusses the results. Finally, Section 7 concludes the paper.

## 2. Proposed Framework and Model

### 2.1. Problem Statement

DERs have the potential to reduce peak demand and enhance energy self-sufficiency. While the integration of renewable-based distributed generation can improve the secure supply of energy for prosumers, it also introduces power grid fluctuations and uncertainties, posing challenges for system operators. P2P energy trading, combined with the flexible operation of DERs, can mitigate sharp demand peaks and deviations. Furthermore, it incentivizes prosumers to contribute their flexible resources, ultimately reducing the overall social cost of energy systems. Fig. 1 illustrate the proposed distributed DRJCC optimization framework.

### 2.2. Energy Community

The high integration of distributed renewable energy resources among prosumers enables P2P energy trading, fostering direct energy sharing and coordinated use of distributed flexibility resources, which benefits the broader power system. However, uncertainties such as weather variability and diverse consumption patterns affect both energy generation and consumption. Ensuring energy balance at all levels—individual prosumers, energy communities, and the wider grid—is critical for achieving flexibility. A key objective remains maintaining energy sufficiency, which can be expressed as

$$\sum_{n \in \mathbb{N}} p_n + q_n + p_n^g + p_n^b + P_n^e \geq \sum_{n \in \mathbb{N}} p_n^l + p_n^s \quad (1)$$

### 2.3. Prosumers Model

Prosumers, equipped with diverse energy resources and managed by Home Energy Management Systems (HEMS), play a pivotal role in distributed energy systems. By actively managing their energy levels and participating in grid operations, prosumers optimize resource scheduling, maximize benefits, and enhance grid stability. Prosumer decision-making considers several factors including their resources, preferences, energy prices, and neighbor's status. Flexible resources such as battery storage, shiftable loads, as well as peer-to-peer energy trading can mitigate grid imbalances caused by fluctuating renewable generation and demand. Even without direct energy imbalances, prosumers can contribute to grid stability and enhance system efficiency by sharing their resources with their peers.

The prosumer's objectives are often framed in terms of their energy purchase $C_n^{da}$ and operational cost $C_n^{op}$, represented mathematically as follows

$$\min_{z \in \mathbb{Z}} J_n(z, \xi) = C_n^{da}(z, \xi) + C_n^{op}(z, \xi) \quad (2)$$

Here, $z$ denotes the set of all decision variables, and $\xi$ represent uncertain model parameters such as PV generation and inflexible loads. Prosumers can generate profit by selling energy to the grid or to their energy partners (the day before), as defined by equation

$$C_n^{da} = c_p^{\top} p_n + \sum_{m \in \mathbb{N}_n} c_{nm}^{\top} p_{nm} \quad (3)$$

Here, the coefficients $c_p$ and $c_{nm}$ are positive and represent the cost of energy exchange with the grid and the other prosumers. The prosumer also incurs costs for battery storage system degradation, the dis-utility caused from shifting usage of flexible loads, and real-time energy purchases, as shown in

$$C_n^{op} = c_q^{\top} q_n + \gamma_n^b \left\| p_n^b \right\|^2 + \gamma_n^s \left\| S_n \right\|^2 \quad (4)$$

In this context, $\gamma_n^b$, $\gamma_n^s$ are positive coefficients. $S_n$ is a variable representing the energy shift state of the flexible loads between desired power $p_n^{s,r}$ and actual values $p_n^s$. For more details, please refer to [18]. It is defined as follows

$$S_{n,t+1} = S_{n,t} + \Delta t(p_{n,t}^s - p_{n,t}^{s,r}), \ \forall t \in T \quad (5)$$

where $\Delta t$ is the sampling time interval. The state-of-charge (SoC) of battery storage systems is defined as

$$E_{n,t+1} = E_{n,t} + \eta p_{n,t}^b, \ \forall t \in T \quad (6)$$

where $\eta$ is the charging/discharging efficiency. To ensure energy balance in the network during P2P energy trading, the following reciprocity constraints should be considered by prosumers

$$p_{nm}^e + p_{mn}^e = 0, \quad m \in \mathrm{N}_n \quad (7)$$

Some other constraints on decision variables and system parameters, due to technical or contractual limitations, are as follows

$$\underline{p}_n \le p_n \le \overline{p}_n \quad (8a)$$

$$\underline{q}_n \le q_n \le \overline{q}_n \quad (8b)$$

$$\underline{p}_{nm}^e \le p_{nm}^e \le \overline{p}_{nm}^e, \quad m \in \mathrm{N}_n \quad (8c)$$

$$\underline{p}_n^b \le p_n^b \le \overline{p}_n^b \quad (8d)$$

$$0 \le p_n^s \le \overline{p}_n^s \quad (8e)$$

$$\underline{E}_n \le E_n \le \overline{E}_n \quad (8f)$$

$$\underline{S}_n \le S_n \le \overline{S}_n \quad (8g)$$

## 3. Problem Formulation

This section addresses uncertainty within both the objectives and constraints of P2P energy trading and flexibility models by exploring distributionally robust optimization and joint chance-constrained techniques.

### 3.1. Distributionally Robust Joint Chance-Constraints

Prosumers often lack an exact probabilistic model P for power generation or energy consumption but have access to a finite set of historical data points $\{\xi_i, i \le I\}$, sampled independently from probability distribution function P. A common practice in such scenarios is to approximate P using the empirical uniform distribution $P_I$, and solve the deterministic sample average approximation (SAA) of the problem. However, this approach often leads to biased decisions with poor out-of-sample performance.

To address this limitation, a data-driven Wasserstein-based DRO framework can be employed. This method hedges against all distributions within a specified Wasserstein distance from $P_I$, represented as ambiguity set $\mathbb{P}$, offering improved robustness and efficiency. Thus, the problem (2) can be expressed as

$$\min_{z \in Z} \sup_{P \in \mathbb{P}} \mathrm{E}_P[J(z, \xi)] \quad (9)$$

where $\mathrm{E}_P[J(z,\xi)]$ is the expected value of $J(z,\xi)$ under the distribution $P$.

On the other hand, some prosumers may fail to meet their P2P scheduled transactions due to uncertain generation and consumption, leading to the risk of cascading constraints violations. Therefore, the DRO problem needs to consider the mis-estimation risks. The joint chance constraints ensure that the maximum simultaneous local power imbalance of prosumers does not exceed required energy of energy community with a probability at least $1 - \varepsilon$, where $\varepsilon \in (0,1)$ is specified risk tolerance. This linear chance-constraints can be given as

$$\inf_{P \in \mathbb{P}} \mathrm{P}\{\xi : a_n\xi \le p_n, \forall n \in N\} \ge 1 - \varepsilon \quad (10)$$

This also enhances energy trading consistency (ETC); further details can be found in [18].

### 3.2. Data-driven Ambiguity Set

In DRO, the ambiguity set plays a crucial role in defining the set of plausible probability distributions for the uncertain parameters. Consider an uncertain parameter vector $\xi$ whose distribution is unknown but is assumed to belong to a set of distributions $\mathbb{P}$. To evaluate the expected value given in (9) or chance-constraints (10), its probability distribution is required. However, in practical applications, the true distribution $P$ is often unknown, and only a set of historical samples $\xi = \{\xi_1, \ldots, \xi_I\}$ is available. In this paper, the Wasserstein metric is employed to construct an ambiguous set $\mathbb{P}$, as it offers desirable properties such as out-of-sample performance guarantees, asymptotic guarantees, and analytical tractability, allowing for a tractable reformulation of the problem [15]. Given a set of historical samples, an empirical distribution $\hat{P}^I = I^{-1} \sum_{i=1}^I \delta_{\hat{\xi}_i}$ can be used to estimate $P$, where $\delta_{\hat{\xi}_i}$ denotes the Dirac measure at point $\hat{\xi}_i$, and $I$ denotes the number of samples. Generally, the Wasserstein metric quantifies the distance between the empirical and true distributions, defined as

$$W(P, \hat{P}^I) = \min_\Pi \left\{ \int_{R^T \times R^T} \|\xi - \bar{\xi}\| \Pi(d\xi, d\bar{\xi}) \right\}$$

Here, $\Pi$ is a joint distribution on $\mathbb{R}^T \times \mathbb{R}^T$ with marginal distributions $P$ and $\hat{P}^I$. Therefore, the ambiguous set can be constructed as

$$\mathrm{P}_I = \left\{ W(P, \hat{P}^I) \le \rho(I) \right\} \quad (11)$$

where $\rho(I)$ is the radius of the ambiguous set centered at $\hat{P}^I$ [16].

### 3.3. Integrated Formulation of P2P Energy Trading and Flexibility Provision under Uncertainty

Prosumers seek to engage in local energy trading to leverage the benefits of renewable, low-cost energy. This mutually advantageous arrangement fosters collaboration between buyers and sellers to optimize flexible resources for enhanced local energy trading. Additionally, transactions that support network operations, such as efficient load shifting to address peak demand, may receive backing from network operators, helping to overcome regulatory challenges. The coordination of P2P energy trading and distributed flexibility provision is modelled as a two-stage distributionally robust optimization problem with chance constraints, as described below.

$$\min \sum_{n \in \mathrm{N}} \left\{ c_p^T p_n + \sum_{m \in \mathrm{N}_n} c_{nm}^T p_{nm}^e + \gamma_n^b \left\| p_n^b \right\|^2 + \gamma_n^s \left\| S_n \right\|^2 \right.$$

$$\left. - \max_{P \in \mathbb{P}} \mathrm{E}_P[c_q^T q_n(\xi_n)] \right\} \quad (12a)$$

$$\text{over } \left\{ p_n, \left\{ p_{nm}^e, \forall m \right\}, p_n^b, p_n^s, q_n(\xi_n), \forall n \in \mathrm{N} \right\}$$

subject to

$$P\left[\begin{matrix}\sum_{n\in N}p_n^l(\xi_n)-p_n^g(\xi_n)\leq\\\sum_{n\in N}p_n+q_n+p_n^b-p_n^s+P_n^e\end{matrix}\right]\geq 1-\varepsilon \quad (12b)$$

$$S_{n,t+1}=S_{n,t}+\Delta t(p_{n,t}^s-p_{n,t}^{s,r}), \quad \forall n\in N, t\in T$$

$$E_{n,t+1}=E_{n,t}+\eta\,p_{n,t}^b, \quad \forall n\in N, t\in T$$

$$p_{nm}^e+p_{mn}^e=0, \qquad \forall n\in N, m\in N_n$$

$$\underline{p}_{nm}^e\leq p_{nm}^e\leq\overline{p}_{nm}, \quad \forall n\in N, m\in N_n$$

$$\underline{p}_n\leq p_n\leq\overline{p}_n, \quad \forall n\in N$$

$$\underline{q}_n\leq q_n\leq\overline{q}_n, \quad \forall n\in N$$

$$\underline{p}_n^b\leq p_n^b\leq\overline{p}_n^b, \quad \forall n\in N$$

$$0\leq p_n^s\leq\overline{p}_n^s, \quad \forall n\in N$$

$$\underline{E}_n\leq E_n\leq\overline{E}_n, \quad \forall n\in N$$

$$\underline{S}_n\leq S_n\leq\overline{S}_n, \quad \forall n\in N$$
$$\left.\right\} \quad (12c)$$

The objective function (12a) comprises two components: the first focuses on planning DA energy purchases, P2P energy trading, and strategic adjustments of flexible resource set-points based on forecasts (i.e., *here-and-now*), while the second addresses real-time power balance adjustments via grid energy purchases (i.e., *wait-and-see*). The energy sufficiency constraint (12b) for prosumers depends on renewable generation and uncertain consumption, and is formulated as a chance constraint to account for potential disruptions that may cause cascading effects.

**Remark1.** The problem (12) still poses significant computational challenges due to the presence of an expected value objective and chance constraints. Both elements require accurate knowledge of the underlying probability distribution P, making the problem intractable in its original form. The following section outlines the assumptions used to reformulate this problem into a tractable form.

## 4. Tractability and Scalability Methodology
In this section, we focus on developing tractable reformulations that enable scalable and distributed implementations.

### 4.1. Reformulation of the Objective Function
To reduce the complexity of problem (12) arising from infinite-dimensional optimization, a common approach [16] is to approximate the functional recourse decisions $q_n(\xi_n)$ in the objective function (12a) using linear decision rules of the form $Q_n\xi_n$, where $Q_n$ is a finite-dimensional coefficient matrix. Therefore, the objective function (12a) can be rewritten as

$$\min\ c_p^T p_n+\sum_{m\in N_n}c_{nm}^T p_{nm}^e$$
$$+\gamma_n^b\left\|p_n^b\right\|^2+\gamma_n^s\left\|S_n\right\|^2-\max_{P\in P}E_P[c_q^T Q_n\xi] \quad (13a)$$

Now, the worst-case expectations of the affine loss function accept a conic program reformulation holding the strong duality and assuming special form uncertainty set. This is stated in the following proposition without a detailed proof; for further information, the reader is referred to [15].

**Proposition 1.** [15] Suppose that the uncertainty set is a polytope, that is, $\Xi=\{\xi\in\mathbb{R}^m:\mathbf{C}\xi\leq d\}$ where $\mathbf{C}$ is a matrix and $d$ a vector of appropriate dimensions. Moreover, consider the linear function $q_n(\xi_n)=\mathbf{Q}_n\xi_n$. Then, the worst-case expectation in (12a) evaluates to

$$\max_{P\in P}E_P[c_q^T Q_n\xi_n]=$$

$$\begin{cases}\min_{\lambda,s,\gamma}\lambda\rho+\dfrac{1}{I}\sum_{i=1}^I s_i\\[2mm]\text{subject to}\quad c_q^T\mathbf{Q}_n\hat{\xi}_i+\gamma_i^T\left(d-\mathbf{C}\hat{\xi}_i\right)\leq s_i, \ \forall i\leq I\\[2mm]\qquad\qquad\left\|\mathbf{C}^T\gamma_i-\mathbf{Q}_n c_q\right\|_*\leq\lambda, \qquad\forall i\leq I\\[2mm]\qquad\qquad\gamma_i\in\mathbb{R}_+, \qquad\qquad\qquad\forall i\leq I\\[2mm]\qquad\qquad\lambda\in\mathbb{R}_+, s\in\mathbb{R}^I \qquad\qquad\forall i\leq I\end{cases} \quad (14a)$$

Here, $\lambda$ denotes auxiliary penalty variable. The corresponding auxiliary epigraph and dual variables are represented by $s_i$ and $\gamma_i, \forall i\leq I$, respectively. $\|.\|_*$ stands for the dual norm of $\|.\|$.

### 4.2. Reformulation of Chance-Constraints
The joint chance constraints (12b) are nonconvex and difficult to solve. Given a set of individual violation tolerances $\varepsilon_n\geq 0$, with $\sum_{n=1}^N\varepsilon_n=\varepsilon$, one can exploit Bonferroni's inequality [7] to split the original joint chance constraint (12b), up into a family of $N$ simpler but more conservative individual chance constraints as follows

$$\min_{P\in P}P\left[p_n^l(\xi_n)-p_n^g(\xi_n)\leq p_n+q_n+p_n^b-p_n^s+P_n^e\right]\geq 1-\varepsilon_n, \forall n \quad (13b)$$

These constraints involve $N$ linear chance-constrained inequalities. For ease of exposition, the left hand side (i.e., the difference of inflexible loads and power generation) can be modeled as $\mathbf{D}_n(\mu_n+\xi_n)$, where $\mathbf{D}_n\in\mathbb{R}^{T\times T}$ is the diagonal matrix of nominal uncertain power of prosumer $n$, $\mu_n\in\mathbb{R}^T$ is the relative predicted uncertain power at first stage, and $\xi_n\in\mathbb{R}^T$ is the uncertain deviation from $\mu_n$, which is revealed at second stage. Optimizing the reformulated individual chance constraints remains challenging, leading us to approximate these worst-case chance constraints with worst-case CVaR constraints, which can be further reformulated into a tractable form, as given by

$$\min_{P\in P}P-\mathrm{CVaR}_{\varepsilon_n}\left[\mathbf{D}_{n,t}(\mu_{n,t}+\xi_{n,t})-(p_{n,t}+q_{n,t}+p_{n,t}^b-p_{n,t}^s+P_{n,t}^e)\right]\leq 0 \quad (15)$$

This approximation is exact for a sufficiently large number of samples, i.e., when $\varepsilon_n\leq I^{-1}$ for all $n\leq N$ [15]. The following result states the reformulation of the CVaR constraint.

**Proposition 2.** The CVaR constraint (15) is equivalent to the following set of constraints.

$$\left\{\min_{P\in P}P-\mathrm{CVaR}_{\varepsilon_n}[\mathbf{D}_n(\mu_n+\xi_n)-(p_n+q_n+p_n^b-p_n^s+P_n^e)]\leq 0\right\}$$

$$=\begin{cases}\lambda_n^c\rho+\dfrac{1}{I}\sum_{i=1}^I s_{n,i}^c\leq 0, \qquad\qquad\forall t\in T\\[2mm]\tau_{n,i}\leq s_{n,i}^c, \qquad\qquad\qquad\forall i\leq I, t\in T\\[2mm]\mathbf{D}_n\xi_{n,i}+\mathbf{D}_n\mu_n-(p_n+p_n^b-\mathbf{D}_n\mu_n-p_n^s+P_n^e)\\[1mm]\quad+(\varepsilon_n-1)\tau_n+\varepsilon_n\gamma_{n,i}^{c\,T}\left(d-\mathbf{C}\hat{\xi}_i\right)\leq s_{n,i}^c, \qquad\forall i\leq I, t\in T\\[2mm]\left\|\varepsilon_n\mathbf{C}^T\gamma_{n,i}^c-\mathbf{D}_n\right\|_*\leq\varepsilon_n\lambda_n, \qquad\forall i\leq I, t\in T\\[2mm]\gamma_{n,i}^c\in\mathbb{R}_+^T, \qquad\qquad\qquad\forall i\leq I, t\in T\\[2mm]\tau\in\mathbb{R}^{N\times T}, \lambda^c\in\mathbb{R}^{N\times T}, s^c\in\mathbb{R}^{N\times I\times T}\end{cases}$$

(14b)

Here, $\lambda_n^c$, $s_{n,i}^c$ and $\gamma_{n,i}^c$, $\forall i \leq I$, denote auxiliary penalty, epigraph and dual variables, respectively.

*Proof.* Refer to Appendix.

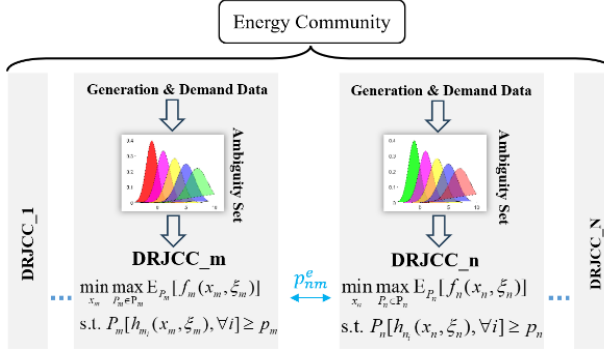

Fig. 2 The proposed scalable and tractable distributed optimization framework

---

**Table 2. Overview of the Distributed Algorithm**

**Initialize** Set $z_n^1$, $\{p_{nm}^1, \forall n \in N_n\}$, $\{\lambda_{nm}^1, \forall n \in N_n\}$, $k = 1$.

**Repeat** the following steps for each prosumer $n$:

- **Update** $k = k + 1$

- **Primal Update:** Solve the optimization problem for the primal variables.

$$z_n^{k+1} = \arg\min_{z_n \in Z_n} L_n(z_n, \hat{p}_{nm}, \lambda_{nm})$$

- **Auxiliary Update:** Update using the closed-form solution.

$$\hat{p}_{nm}^{k+1} = -\hat{p}_{mn}^{k+1} = \frac{p_{nm}^{k+1} - p_{mn}^{k+1}}{2} + \frac{\lambda_{nm}^k - \lambda_{nm}^k}{2\sigma}$$

- **Dual Update:** Update the dual variables using the dual feasibility condition.

$$\lambda_{nm}^{k+1} = \lambda_{nm}^k + \sigma\left[p_{nm}^{k+1} - \hat{p}_{nm}^{k+1}\right]$$

**Until** $(\Delta e_p^{k+1} \leq \varepsilon_p$ and $\Delta e_D^{k+1} \leq \varepsilon_D)$

The residuals are defined as

$$\Delta e_p^{k+1} = \sum_{m \in N_n} \left\| p_{nm}^{k+1} - \hat{p}_{nm}^{k+1} \right\|$$

$$\Delta e_D^{k+1} = \sum_{m \in N_n} \left\| \lambda_{nm}^{k+1} - \lambda_{nm}^k \right\|$$

---

### 4.3. Distributed Implementation based on ADMM

The reformulated problem (14) along with constraints (12c) can be solved using a distributed approach by decoupling the complicating reciprocity constraints (7). However, the problem results in a multi-block alternating direction method of multipliers (ADMM), which, as highlighted in the literature (e.g., [20]), is known for its slow convergence and, the lack of guaranteed convergence. To mitigate these challenges, we introduce auxiliary variables, following well-established strategies (e.g., [8]), which allow the problem to be reformulated into a two-block ADMM. This transformation not only improves the convergence speed but also enhances the stability of the algorithm by expressing the reciprocity constraints (7) as

$$\hat{p}_{nm} = p_{nm} : \lambda_{nm}, \quad \forall n \in N, m \in N_n \quad (16a)$$

$$\hat{p}_{nm} + \hat{p}_{mn} = 0, \quad \forall n \in N, m \in N_n \quad (16b)$$

The augmented Lagrangian of problem (14) is formulated as

$$L(z_n, p_{nm}, \lambda_{nm}) =$$
$$\sum_{n \in N} J_n(z_n, p_{nm}) + \lambda_{nm}^T(p_{nm} - \hat{p}_{nm}) + \sigma/2 \left\| p_{nm} - \hat{p}_{nm} \right\|^2$$

where $\lambda_{nm}$ denotes the Lagrangian multipliers of (16b); $\sigma$ is the penalty factor in the ADMM algorithm. The ADMM algorithm for solving the coordination problem are given in Table 2. The closed-form solution of auxiliary variables update is derived by solving Karush-Kuhn-Tucker (KKT) conditions of optimization of augmented Lagrangian constrained by (7b) over auxiliary variables. The derivation is provided in Appendix. For further details, the reader is referred to [18]. Fig. 2 illustrates the proposed scalable and tractable distributed coordination framework.

## 5. Numerical Implementation

### 5.1. Case Study

To evaluate the proposed model's performance, we conducted an extensive case study on a distributed energy system comprising ten prosumers in central Iran. The system incorporates two main sources of uncertainty: PV generation and inflexible loads. Their nominal values, along with resource capacity limits and trading specifications, are detailed in Table 3 [18]. The electricity pricing structure includes both DA and RT markets, following peak/off-peak periods and time-of-use (TOU) rates, as illustrated in Fig. 3. The study utilizes hourly electricity consumption and solar generation data collected during Summer 2019-2021 [21], with comprehensive data analysis presented in the following section. The proposed distributionally robust optimization model was implemented in MATLAB using the MOSEK solver. All computations were performed on a workstation equipped with an Intel Core i7-11700 CPU (2.5 GHz, 8 cores) and 16 GB RAM.

Table 3. Typical Data and Model parameters

| Parameters | |
|---|---|
| $\rho = [0 \ 1]$ | $\underline{p}_n = -\overline{p}_n \in [40 \ 60]$ kW |
| $\sigma = 0.1$, $\varepsilon = 10^{-6}$ | $\underline{p}_{nm}^e = -\overline{p}_{nm}^e = 5$ kW, $\forall m \in N_n$ |
| $C = 0, d = 0, \Delta t = 1$ | $\underline{p}_n^b = -\overline{p}_n^b \in [10 \ 40]$ kW |
| $D = [-10 \ 10]$ kW | $\overline{p}_n^s \in [10 \ 40]$ kW |
| $E \in [-200 \ 200]$ kWh | $c_{nm} \in [4 \ 16]¢$ |
| $S \in [0 \ 200]$ kWh | $c_p \in [6 \ 12]¢$, $c_q \in [6 \ 25]¢$ |

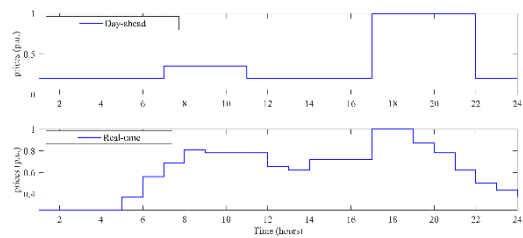

Fig. 3. Normalized DA (top) and RT (bottom) energy prices.

### 5.2. Data Analysis Results

To assess the dataset's properties, we analysed data from 100 prosumers and applied k-means clustering to identify distinct consumption patterns. The optimal number of clusters (k) was determined using the Elbow method, which analyses the within-cluster sum of squares (WCSS) as a function of k [19]. The Elbow method suggested an optimal k of 4, with a corresponding WCSS of 600. Based on their observed characteristics, these four clusters were subsequently labelled as "Residential", "Commercial", "Industrial", and "Public". To further characterize the identified consumer clusters, a multi-dimensional analysis was conducted. The following key load profile indicators, as commonly employed in load profiling studies [17], were considered:

- *Night Consumption Ratio* (NCR) (10PM-6AM): Reflects the proportion of energy consumed during off-peak hours.

- *Business Hours Ratio* (BHR) (8AM-6PM): Represents the share of energy consumed during typical business hours.

- *Load Factor* (LF): Indicates the average load as a fraction of the peak load, reflecting the consistency of energy demand.

- *Consumption Variance* (CV): Measures the variability of energy consumption over time.

- *Peak Hour*: Identifies the hour of maximum energy consumption.

- *Peak-to-Average Power Ratio* (PAPR): Quantifies the ratio between peak demand and average demand.

- *Daily Total Consumption* (DTC): Represents the overall energy consumption per day.

The identified features enabled the differentiation of four distinct consumption patterns, as illustrated in Fig. 5. This analysis provides valuable insights into the energy profiles of these patterns, facilitating informed decision-making in demand-side management, energy market design, and resource allocation optimization for enhanced grid stability.

The average 24-hour consumption profiles for each pattern are presented in Fig. 6, while Fig. 7 depicts the distribution of hourly peak demand. A comprehensive summary of the statistical characteristics is provided in Table 4.

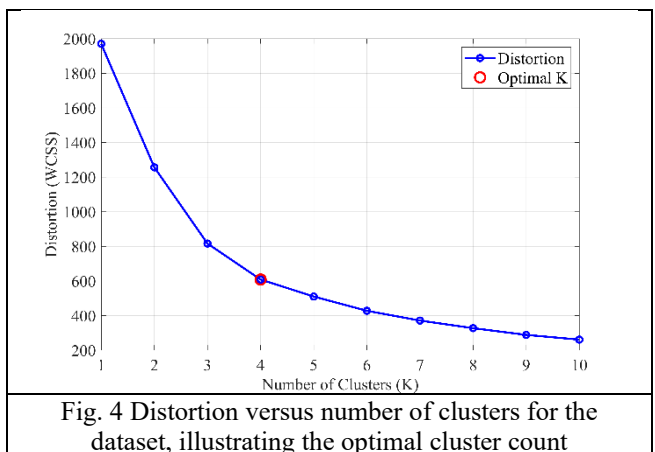

Fig. 4 Distortion versus number of clusters for the dataset, illustrating the optimal cluster count

Table 4 Statistical charactristics (Averaged)

| Patterns | Average Daily consumption (kWh) | Peak-to-average ratio | Night consumption ratio | Business hours ratio | Load factor | Most Common Peak hours |
|---|---|---|---|---|---|---|
| Residential | 6.99 | 3.32 | 0.49 | 1.06 | 0.31 | 18 |
| Commercial | 16.61 | 2.43 | 0.25 | 1.73 | 0.42 | 11 |
| Industrial | 28.24 | 1.76 | 0.53 | 1.36 | 0.57 | 17 |
| Public | 7.52 | 2.51 | 0.32 | 1.37 | 0.40 | 18 |

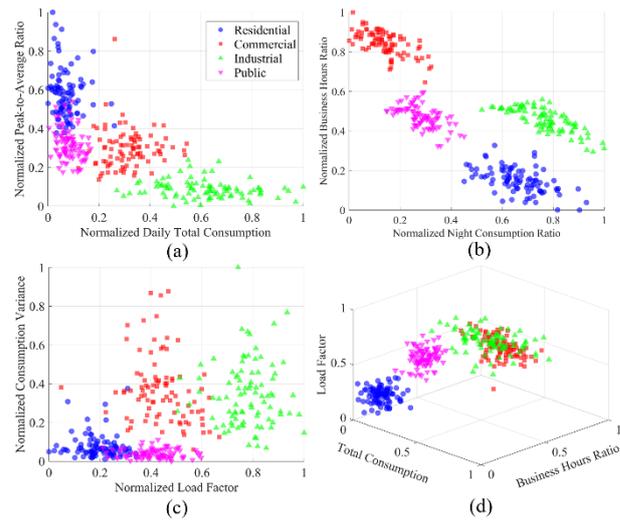

Fig. 5. Feature space exploration of the dataset; (a) Normalized PAPR vs. DTC, (b) Normalized BHR vs. NCR, (c) Normalized CV vs. LF, and (d) Normalized LF vs. TC-BHR

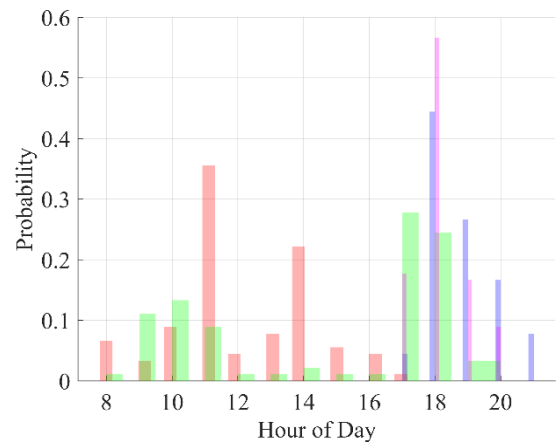

Fig. 6 Hourly peak demand distribution.

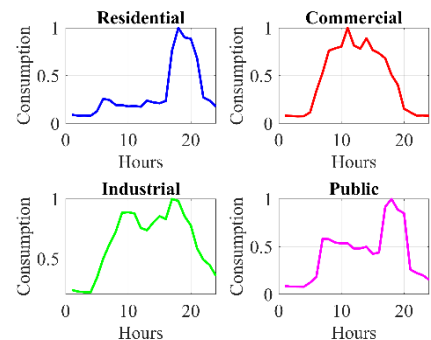

Fig. 7 Average consumption power of various patterns

### 5.3. Implementation Results

Figures 8 and 9 show the energy profiles of two representative prosumers (residential and commercial) over a 24-hour period, comparing the baseline scenario with the proposed approach. In the baseline scenario, no P2P energy trading occurs, and flexible resources, such as BES and SLs, are adjusted independently to maintain power balance for each prosumer in a stand-alone optimization framework. The operation of flexible resources is suboptimal due to limited options and trade-offs. In contrast, the proposed strategy enables P2P energy trading, allowing both explicit energy exchange and implicit flexibility sharing.

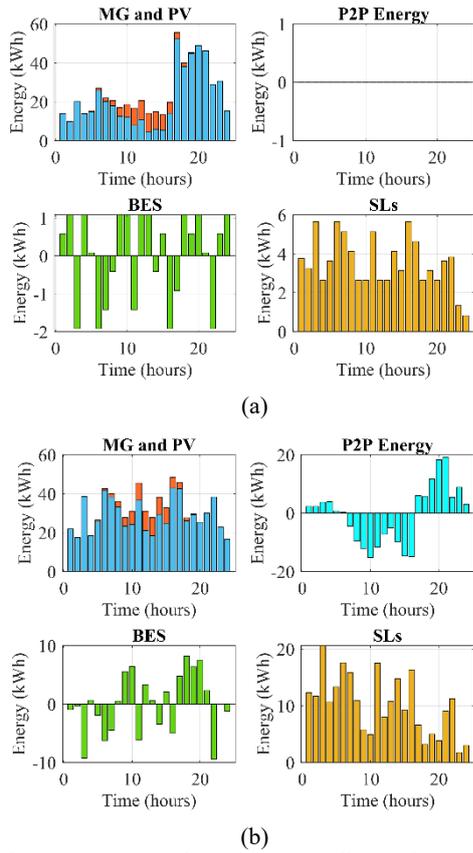

(a)

(b)

Fig. 8. Energy generation, storage, trading and consumption of a representative residential prosumer, (a) Baseline condition, (b) The proposed strategy

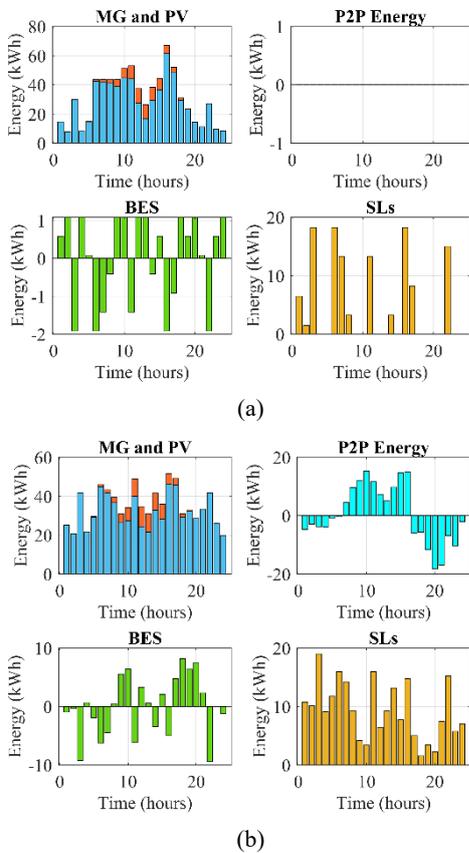

(a)

(b)

Fig. 9. Energy generation, storage, trading and consumption of a representative commercial prosumer, (a) Baseline conditions, (b) The proposed strategy

When P2P energy trading occurs between prosumers, flexible resources are efficiently utilized in collaboration with other prosumers having different patterns and operations. This is particularly evident in the coordination of SLs and BES, which are optimized to manage PV generation during the midday peak.

Table 5 highlights a 28% reduction in total cost with the proposed strategy, driven by the efficient and collaborative use of flexible resources and effective management of uncertainties.

Additionally, the strategy achieves a 31% reduction in peak demand, primarily due to the coordinated use of flexible resources through P2P energy trading. Figures 10(a) and 10(b) show the energy demand of prosumers over a 24-hour period, before and after the implementation of the proposed approach, illustrating a substantial decrease in peak demand.

Table 5 Comparison of prosumers total cost ($)

|  | Baseline Conditions | Proposed Strategy | Relative Reduction (%) |
|---|---|---|---|
| Total cost ($) | 744.21 | 533.79 | 28.3 |
| PAR (%) | 3.8 | 2.7 | 31.7 |

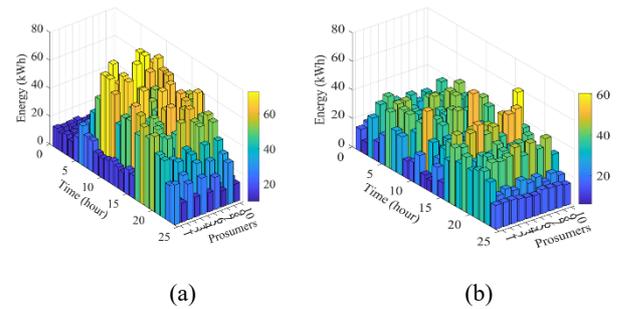

(a)                                    (b)

Fig. 10. Energy demand and peak load profiles of prosumers: (a) under baseline conditions, and (b) with the implementation of the proposed strategy.

## 6. Discussions

### 6.1. Out-of-Sample Analysis

A key advantage of data-driven DRO is its ability to handle uncertainty without relying on specific probability distributions, making it effective in dynamic and ambiguous environments. We evaluated our framework by training the model on 2019 data and testing it on 100 samples from 2020-2021, assessing out-of-sample costs and constraint violations with the sample average approximation (SAA) method. Figure 7 illustrates the impact of the Wasserstein radius ρ on the cost-violation trade-off. Smaller ρ reduces costs but increases violation risk, while larger ρ results in conservative, more expensive strategies. At large ρ values, the empirical violation probabilities remain below 1%, ensuring robust constraint satisfaction.

Table 6 presents the average out-of-sample total cost and within-cluster costs, showing a decrease in costs as ρ decreases from 0.2 to 0.01, followed by a slight increase at ρ = 0.001, suggesting over-conservatism. However, inter-cluster analysis reveals that the optimal ρ radius of the ambiguity set varies across clusters and scenarios, with cluster-specific adjustments offering potential for further cost reductions.

$$\hat{C}(\rho) = \sum_{n \in N} \left\{ c_p^T \hat{p}_n(\rho) + \sum_{m \in \mathcal{N}_n} c_{nm}^T \hat{p}_{nm}^e(\rho) \right.$$

$$\left. + \gamma_n^b \left\| p_n^b(\rho) \right\|^2 + \gamma_n^s \left\| S_n(\rho) \right\|^2 - \max_{P \in P} \mathrm{E}_P[c_q^T Q_n(\rho)\xi_n] \right\}$$

$$\hat{V}(\rho) = \frac{1}{I} \sum_{n \in N} \sum_{i=1}^{I} \mathbf{1}_{(C_n(\mu_n + \xi_n) \le \hat{p}_n(\rho) + \hat{q}_n(\rho) + \hat{p}_n^b(\rho) - \hat{p}_{nm}^s(\rho) + \hat{p}_n^{ g}(\rho))}$$

Table 6. Average out-of-sample costs

| $\rho$ | 0.2 | 0.1 | 0.03 | 0.01 | 0.001 |
|---|---|---|---|---|---|
| $C^{res}(\rho)$ | 72.01 | 68.73 | 66.37 | **60.31** | 65.24 |
| $C^{com}(\rho)$ | 177.72 | 165.36 | **149.48** | 151.52 | 155.96 |
| $C^{ind}(\rho)$ | 297.64 | 279.49 | **254.46** | 257.14 | 265.11 |
| $C^{pub}(\rho)$ | 78.15 | 73.68 | 70.75 | **65.16** | 71.84 |
| $\hat{C}(\rho)$ | 624 | 587 | 541 | 534 | 557 |

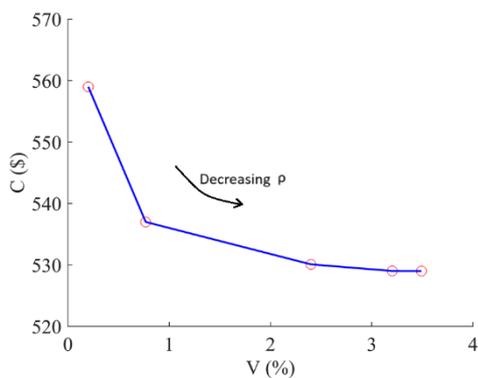

Fig. 11 Out-of-sample costs and violation probabilities

### 6.2. Distributed algorithm results

The ADMM algorithm was used for distributed implementation, demonstrating steady convergence behavior as seen in the residual plots in Fig. 12. The cost optimization results show a consistent decrease in the objective function, indicating successful minimization of overall costs. The algorithm achieves an acceptable tolerance within fewer than 30 iterations, with each iteration taking approximately 1.89 seconds.

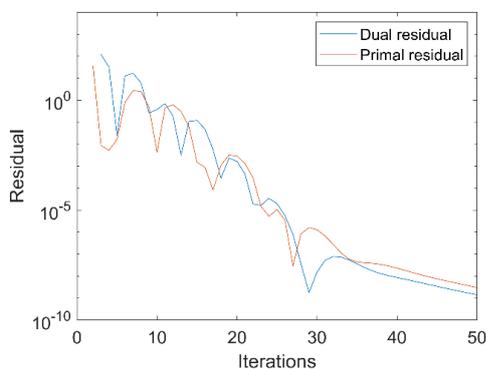

Fig. 12 Primal and dual residual convergence

### 6.3. Cross-Market Adaptability of the Decentralized Framework

The proposed model demonstrates strong adaptability across diverse consumption patterns and market environments, as evidenced by its ability to achieve 28% cost reduction in energy trading under varying pricing schemes. Leveraging the DRO approach, the model facilitates seamless adjustments to energy pricing mechanisms and collective incentive structures, maintaining up to 31% peak demand reduction across scenarios. These results highlight its broad applicability, particularly in regions with high renewable energy penetration. Furthermore, the decentralized coordination enabled by the ADMM algorithm effectively addresses scalability challenges. The algorithm demonstrates robust performance by handling 10 distinct consumption patterns among prosumers with 99.9% convergence accuracy in under 20 seconds. This efficiency, combined with the model's adaptability through parameter tuning and flexible resource utilization, ensures its suitability for both developed and developing regions with diverse grid configurations. These features underscore the model's practical applicability in real-world energy systems.

## 7. Conclusion

This study addresses critical challenges in power grids under uncertainty, particularly the growing peak demand and capacity gap. The proposed framework incentivizes prosumers by reducing costs and enhancing the resilience of energy communities, while simultaneously providing grid services. To achieve this, a novel two-stage, multi-period distributionally robust optimization framework with joint chance constraints is introduced to manage prosumer operations and energy sharing, mitigating peak load imbalances under uncertainty. Several tractable reformulations, based on the Wasserstein metric and CVaR, are employed to enable distributed implementation. Additionally, machine learning techniques are applied to analyse the feature space of interacting agents, identifying potential collaborations and guiding the adjustment of ambiguity set parameters. A case study demonstrated a 28% reduction in power imbalances and a 31% decrease in prosumer costs. The model's integration of flexible loads and energy storage enhances system stability, presenting a promising solution for grid resilience and efficiency, with applicability across various regulatory frameworks and market conditions.

## 8. References

[۱] ح. بیات، ف. اسدی، "بررسی ابعاد امنیت تأمین برق در اوج مصرف برق تابستان۱۴۰۳"، مرکز پژوهش های اتاق ایران، مرداد ۱۴۰۳.

[۲] ا. محمدی، ن. رضایی، م. غلامی، "چارچوب بهینه تجارت انرژی نظیر به نظیر برای ریزشبکه های DC به هم پیوسته با در نظر گرفتن محدودیت تلفات توان"، مجله مهندسی برق دانشگاه تبریز، جلد ۵۴، شماره ۳، پاییز ۱۴۰۳.

[۳] م. کجوری نفت چالی، ع. فریدونیان، ح. لسانی، "پیک سایی تطبیقی و داده محور در شبکه هوشمند انرژی الکتریکی با تحلیل داده های زیرساخت اندازه گیری پیشرفته"، مجله مهندسی برق دانشگاه تبریز، جلد ۴۹، شماره ۳، پاییز ۱۳۹۸.

[۴] س. همتی، س. ف. قادری، م. ص. قاضی زاده، "طراحی هاب انرژی پایدار با درنظرگرفتن ریسک با استفاده از الگوریتم تجزیه بندرز"، مجله مهندسی برق دانشگاه تبریز، جلد ۵۰، شماره ۱، بهار ۱۳۹۹.

[5] W. Yuwei, Y. Yang, L. Tang, W. Sun, and B. Li. "A Wasserstein based two-stage distributionally robust optimization model for optimal operation of CCHP

micro-grid under uncertainties", *International Journal of Electrical Power & Energy Systems* vol. 119, 2020.

[6] Z. Guo, P. Pinson, S. Chen, Q. Yang, Z. Yang, "Chance-constrained peer-to-peer joint energy and reserve market considering renewable generation uncertainty" *IEEE transactions on smart grid,* vol. 12, no. 1, pp. 798-809, 2020.

[7] O. Christos, V. A. Nguyen, D. Kuhn, P. Pinson, "Energy and reserve dispatch with distributionally robust joint chance constraints", *Operations Research Letters,* vol. 49, no. 3, pp. 291-299, 2021.

[8] J. Li, M. E. Khodayar, J. Wang, B. Zhou. "Data-driven distributionally robust co-optimization of P2P energy trading and network operation for interconnected microgrids" *IEEE Transactions on Smart Grid,* vol. 12, no. 6, pp. 5172-5184, 2021.

[9] X. Li, C. Li, G. Chen, Z. Y. Dong, "A data-driven joint chance-constrained game for renewable energy aggregators in the local market", *IEEE Transactions on Smart Grid,* vol. 14, no. 2, pp. 1430-1440, 2022.

[10] C. Zhang, H. Liang, Y. Lai. "A distributionally robust energy management of microgrid problem with ambiguous chance constraints and its tractable approximation method" *Renewable Energy Focus,* vol. 48, 100542, 2024.

[11] S. Suthar, N. M. Pindoriya, "Chance-constrained co-optimization of peer-to-peer energy trading and distribution network operations" *Sustainable Energy, Grids and Networks,* vol. 38, 101344, 2024.

[12] J. Zhai, Y. Jiang, M. Zhou, Y. Shi, W. Chen, C. N. Jones, "Data-Driven Joint Distributionally Robust Chance-Constrained Operation for Multiple Integrated Electricity and Heating Systems" *IEEE Transactions on Sustainable Energy,* (Early Access), 2024.

[13] Z. Junjie, Y. Zhao, Y. Li, M. Yan, Y. Peng, Y. Cai, Y. Cao. "Synergistic Operation Framework for the Energy Hub Merging Stochastic Distributionally Robust Chance-Constrained Optimization and Stackelberg Game", *IEEE Transactions on Smart Grid* (Early Access), 2024.

[14] J. Wenhao, T. Ding, Y. Yuan, C. Mu, H. Zhang, S. Wang, Y. He, and X. Sun. "Decentralized distributionally robust chance-constrained operation of integrated electricity and hydrogen transportation networks", *Applied Energy*, vol. 377, 124369, 2025.

[15] P. Mohajerin Esfahani, D. Kuhn, "Data-driven distributionally robust optimization using the Wasserstein metric: Performance guarantees and tractable reformulations", *Mathematical Programming,* vol. 171, no. 1, pp. 115-166, 2018.

[16] X. Peng, P. Jirutitijaroen, C. Singh. "A distributionally robust optimization model for unit commitment considering uncertain wind power generation." *IEEE Transactions on Power Systems,* vol. 32, no. 1, pp. 39-49, 2016.

[17] Y. Wang, Q. Chen, C. Kang, M. Zhang, K. Wang, Y. Zhao, "Load profiling and its application to demand response: A review", *Tsinghua Science and Technology,* vol. 20, no. 2, pp. 117-129, 2015.

[18] A. Noori, B. Tavassoli, A. Fereidunian, "Joint flexibility-risk managed distributed energy trading considering network constraints and uncertainty",

[19] S. M. Miraftabzadeh, C. G. Colombo, M. Longo, F. Foiadelli, "K-means and alternative clustering methods in modern power systems" *IEEE Access,* vol. 11, pp. 119596 - 119633, 2023.

[20] S. Boyd, N. Parikh, E. Chu, B. Peleato, J. Eckstein, "Distributed optimization and statistical learning via the alternating direction method of multipliers" *Foundations and Trends® in Machine learning.* vol. 3, no. 1, pp. 1-122, 2011.

[21] A. Noori, B. Tavassoli, A. Fereidunian, "Incentivizing peer-to-peer energy trading in microgrids" In *2021 29th Iranian Conference on Electrical Engineering (ICEE)*, pp. 323-328. IEEE, 2021.

*Electric Power Systems Research,* vol. 231, 110355 , 2024.

# 9. Appendix

*Proof of* Proposition 2. For ease of exposition, we define

$$a_n \xi + b_n \equiv D_n \xi + (D_n u_n - p_n + q_n + p_n^b - p_n^s + P_n^e)$$

The worst-case CVaR can be rewritten as

$$\max_{p \in P} P - CVaR_{\tau_n} [a_n \xi + b_n]$$

$$= \max_{p \in P} \min_{\beta \in R} \left\{ \beta + \frac{1}{\varepsilon} E_p [a_n \xi + b_n]_+ \right\}$$

$$= \max_{p \in P} \min_{\beta \in R} \left\{ E_p \left[ \max \left\{ \beta, \frac{1}{\varepsilon} (a_n \xi + b_n) + (1 - \frac{1}{\varepsilon}) \beta \right\} \right] \right\}$$

$$= \min_{\beta \in R} \left\{ \max_{p \in P} E_p \left[ \max \left\{ \beta, \frac{1}{\varepsilon} (a_n \xi + b_n) + (1 - \frac{1}{\varepsilon}) \beta \right\} \right] \right\}$$

where $[.]_+ = max(.,0)$. The first equality uses the definition of CVaR at level $\varepsilon$ with probability distribution $P$. The last equality follows from minimax inequality, in case where objective functions are linear and convex constrained within a weakly compact set [7]. Then, the worst-case expectation can be reformulated as a finite convex minimization problem as follows [15]

$$
\begin{cases}
\lambda_{n,i} \rho + \frac{1}{I} \sum_{i=1}^{I} s_{n,i,t} \leq 0, & \forall t \in T \\
\tau_{n,t} \leq s_{n,i,t}, & \forall i \in I, t \in T \\
D_{n,t} \hat{\xi}_{n,i} + D_{n,t} H_{n,t} - (p_{n,t} + p_{n,t}^b - D_{n,t} u_{n,t} - p_{n,t}^s + P_{n,t}^e) \\
\quad + (\varepsilon_{n,i} - 1) \tau_{n,t} + \varepsilon_{n,i} \gamma_{n,i,t}^T (d - C \hat{\xi}_n) \leq s_{n,i,t}, & \forall i \in I, t \in T \\
\| \varepsilon_{n,i} C^T \gamma_{n,i,t} - D_{n,t} \|_* \leq \varepsilon_{n,i} \lambda_{n,t}, & \forall i \in I, t \in T \\
\gamma_{n,i,t} \in R_+, & \forall i \in I, t \in T \\
\tau \in R^{N \times T}, \lambda \in R^{N \times T}, s \in R^{N \times I \times T}
\end{cases}
$$

∎

*Derivation of closed-form solution* (Auxiliary Update). The augmented Lagrangian function for auxiliary variables can be defined as follows

$$L = -\lambda_{nm}^{\ T} \hat{p}_{nm} - \lambda_{mn}^{\ T} \hat{p}_{mn} + \nu^T (\hat{p}_{nm} + \hat{p}_{mn})$$

$$+ \frac{\sigma}{2} \| p_{nm} - \hat{p}_{nm} \|_2^2 + \frac{\sigma}{2} \| p_{mn} - \hat{p}_{mn} \|_2^2$$

From the first-order optimality conditions, we have

$$\frac{\partial L}{\partial \hat{p}_{nm}} = -\lambda_{nm} + \sigma (p_{nm} - \hat{p}_{nm}) + \nu = 0$$

$$\frac{\partial L}{\partial \hat{p}_{mn}} = -\lambda_{mn} + \sigma (p_{mn} - \hat{p}_{mn}) + \nu = 0$$

By substituting the reciprocity constraints, the auxiliary variable update rules are derived as

$$\hat{p}_{nm}^{k+1} = \frac{p_{nm}^{k+1} - p_{mn}^{k+1}}{2} + \frac{\lambda_{nm}^{k+1} - \lambda_{mn}^{k+1}}{2\sigma}$$

∎